\documentclass[11pt,reqno]{amsart}
\usepackage{amsmath}
\usepackage{amsthm}
\usepackage{amssymb}
\usepackage{enumerate}
\usepackage{amscd}
\usepackage{color}
\usepackage{pb-diagram}
\usepackage{graphicx}
\usepackage[all, cmtip]{xy}
\usepackage{soul}
\usepackage{esint}
\usepackage{hyperref}
\hypersetup{pdftex,colorlinks=true,allcolors=blue}
\usepackage{hypcap}
\usepackage[titletoc]{appendix}
\usepackage{comment}

\usepackage{amsrefs}

\theoremstyle{plain}
\newtheorem{lemma}{Lemma}

\newtheorem{theo}[lemma]{Theorem}
\newtheorem{rema}[lemma]{Remark}
\newtheorem{lem}[lemma]{Lemma}

\newtheorem{thm-Intro}{Theorem} 
\newtheorem{cor-Intro}{Corollary} 
\numberwithin{equation}{section}

\newcommand{\Hol}{\textup{Hol}}


\begin{document}
\title[Invariant Distances on Complete Reinhardt Domains]{Comparing the Carath\'eodory Pseudo-Distance and the K\"ahler-Einstein Distance on Complete Reinhardt Domains}

	\author{Gunhee Cho}
\address{Department of Mathematics\\
	University of California, Santa Barbara\\
	Santa Barbara, CA 93106}
\email{gunhee.cho@math.ucsb.edu}

\begin{abstract} 
We show that on a certain class of bounded, complete Reinhardt domains in $\mathbb{C}^n$ that enjoy a lot of symmetries, the Carath\'eodory pseudo-distance and the geodesic distance of the complete K\"ahler-Einstein metric with Ricci curvature $-1$ are different. 
\end{abstract}

\maketitle


\section{Introduction and main results}

In this paper, we compare the Carath\'eodory pseudo--distance and the geodesic distance of the complete K\"ahler--Einstein metric with Ricci curvature $-1$ on certain complete Reinhardt domains in $\mathbb{C}^n,n\geq 2$. Throughout the paper, we will call the geodesic distance of the complete K\"ahler--Einstein metric with Ricci curvature $-1$ by the K\"ahler--Einstein distance.  The Carath\'eodory pseudo--distance $c_{M}$ on a complex manifold $M$ is defined by 
\[c_{M}(x,y):=\sup_{f\in \Hol(M,\mathbb{D})} \rho_{\mathbb{D}}(f(x),f(y)).\]
Here, $\rho_{\mathbb{D}}$ denotes the Poincar\'e distance on the unit disk $\mathbb{D}$ in $\mathbb{C}^1$ and $\Hol(M,\mathbb{D})$ is the collection of all holomorphic functions from $M$ to $\mathbb{D}$. We call $c_{M}$ the \textit{pseudo}-distance instead of the distance because there could be distinct points $x,y \in M$ such that $c_{M}(x,y) = 0$.

The Carath\'eodory--Reiffen (pseudo--) metric and the Carath\'eodory pseudo--distance are objects of interest in the long-standing conjecture on the existence of a bounded, non-constant, holomorphic function on a simply connected, complete K\"ahler manifold whose sectional
curvature is bounded from above by a negative constant. Recently, Wu and Yau proved that on a manifold with negatively pinched sectional curvature, the K\"ahler--Einstein, Bergman, and Kobayashi-Royden metrics exist and are uniformly equivalent to one another \cite{WuDaminYauShingTung20}. Thus, it is natural to compare the three invariant metrics above to Carath\'eodory-Reiffen metric \cite{WuH93,WuYau20}. The existence of the complete K\"ahler--Einstein metric and distance on wide classes of negatively curved complex manifolds that are invariant under biholomorphic mappings have been studied; especially for various classes of pseudoconvex domains in $\mathbb{C}^n$. In recent years a number of results devoted to the study of the relation between the K\"ahler--Einstein metric and other invariant objects were obtained  (see, for example, \cite{WuH93,WuDaminYauShingTung20,GrahamIan87, ChengShiuYuenYauShingTung80,LempertLaszlo81,PflugPeterZwonekWlodzimierz98,NikolovNikolaiPflugPeter08,NikolovNikolaiPflugPeter09,NikolovNikolaiPflugPeter11,NikolovNikolaiPflugPeterThomasPascalJZwonekWlodzimierz08,NikolovNikolaiPflugPeterZwonekWlodzimierz07,YeungSai-Kee09,CatlinDavidW89,HerbortGregor05,KosinskiLukasz14,FuSiqi14,ChoGunhee21,ChoGunheeQianJunqing20,QianJunqing20,WangWei96,BlandJohnS86,BlankBrianEFanDaShanKleinDavidKrantzStevenGMaDaoweiPangMyung-Yull92,GunheeChoYuanYuan20,GunheeCho20,WuYau20} and references therein). The results involved also the smallest invariant metric and the invariant distance among the possible invariant metrics and invariant distances, i.e., the Carath\'eodory-Reiffen metric and Carath\'eodory pseudo--distance (\cite[Proposition 3.1.7]{KobayashiShoshichi98},\cite[Chapter 2]{JarnickiMarekPflugPeter13}).  

As an invariant pseudo--distance, the Carath\'eodory pseudo--distance is less than or equal to the Carath\'eodory inner--distance \cite[Remark 2.7.5]{JarnickiMarekPflugPeter13}, Bergman distance \cite{LookKH57,HahnKyongT76,HahnKyongT77}, and Kobayashi-Royden distance \cite{KobayashiShoshichi98,JarnickiMarekPflugPeter13} and specific examples in which Carath\'eodory pseudo--distance is strictly smaller than each of these distances are known (for example, \cite{VigueJeanPierre83,NikolovNikolaiPflugPeterThomasPascalJZwonekWlodzimierz08,NikolovNikolaiPflugPeterZwonekWlodzimierz07}). As a consequence of the Schwarz--Yau Lemma (Section 2), the Carath\'eodory pseudo--distance is less than or equal to the K\"ahler--Einstein distance but a specific example was not investigated. In this paper, we provide a class of complete Reinhardt domains in $\mathbb{C}^n, n\geq 2$ that distinguish the K\"ahler--Einstein distance from the Carath\'eodory pseudo--distance. 

The difficulty of distinguishing the K\"ahler--Einstein distance from the Carath\'eodory pseudo--distance is that even in the case of strictly pseudoconvex domains in $\mathbb{C}^n$, the concrete formulas of the K\"ahler--Einstein distance and the Carath\'eodory pseudo--distance are unknown. Thus, a concrete comparison of the two distances is non-trivial (see \cite{VenturiniSergio89,LeeJohnMMelroseRichard82} and \cite[Proposition 5.5, Theorem 7.5]{ChengShiuYuenYauShingTung80} for the case of smoothly bounded, strictly pseudoconvex domains). Vigue's Work \cite{VigueJeanPierre83} inspired one approach to compare the Carath\'eodory pseudo--distance and the Carath\'eodory  inner--distance. After a clever choice of one point with origin on the diagonal entries on
\begin{equation}\label{eq:example}
\left\{ (z_1,z_2)\in\mathbb{C}^{2}:\left| z_1 \right|+\left| z_2 \right|<1, \left| z_1 z_2 \right|<\frac{1}{16} \right\},
\end{equation}
Vigue distinguished the Carath\'eodory inner-distance and the Carath\'eodory pseudo--distance. By extending this idea, we could distinguish the Carath\'eodory-distance from the K\"ahler--Einstein distance.  
We say a domain $\Omega$ in $\mathbb{C}^n$ is a complete Reinhardt domain if for any $(z_1,...,z_n)\in \Omega$ and $(\lambda_1,...,\lambda_n)\in \overline{\mathbb{D}}^n$, the point $(\lambda_1 z_1,...,\lambda_n z_n)$ belongs to $\Omega$ \cite[Remark 2.2.1]{JarnickiMarekPflugPeter13}. We say $\Omega \subset \mathbb{C}^n$ is compatible with the symmetry by all permutations if $\Omega$ satisfies the following: $(z_1,\cdots,z_a,\cdots,z_b,\cdots,z_n) \in \Omega$ if and only if $(z_1,\cdots,z_{\phi(a)},\cdots,z_{\phi(b)},\cdots,z_n) \in \Omega$ for any $a,b\in  \left\{1,\cdots,n \right\}$ and any permutation $\phi : \left\{1,\cdots,n \right\} \rightarrow \left\{1,\cdots,n \right\}$.
We denote the ball of radius $R>0$, centered at the origin in $\mathbb{C}^n$ by $B_R$. We denote the complete K\"ahler--Einstein metric with Ricci curvature $-1$ on $\Omega$ by $\omega_{KE}$, and let $d^{{KE}}_{\Omega}$ be the distance induced by $\omega_{KE}$. 

\begin{theo}\label{thm:main}
	Let $\Omega$ be a bounded, complete Reinhardt domain in $\mathbb{C}^n, n\geq 2$. Suppose that $\Omega$ is compatible with the symmetry by all permutations and there exists $R>0$ such that $\Omega$ is contained in $B_R$ and $\sum_{k=1}^{n} p_k=R$ for some $(p_1,\cdots,p_n)\in \Omega$. Then 	
	\begin{equation}\label{eq:6}
	c_{\Omega} \lneq d^{{KE}}_{\Omega}. 
	\end{equation}
Here, \eqref{eq:6} means for any $a,b\in \Omega$, and there exist $p,q\in \Omega$ such that 	
	\begin{equation*}
	c_{\Omega}(a,b) \leq d^{{KE}}_{\Omega}(a,b),
	\end{equation*}
	\begin{equation*}
	c_{\Omega}(p,q) < d^{{KE}}_{\Omega}(p,q).
	\end{equation*}\end{theo}

We apply the theorem of the existence of the complete K\"ahler--Einstein metric of a negative Ricci curvature on a bounded pseudoconvex domain in $\mathbb{C}^n$ as given in the main theorem in \cite{MokNgaimingYauShingTung83}. The hypothesis of Theorem~\ref{thm:main} is satisfied in several examples including ~\eqref{eq:example} complex ellipsoids \cite{ChoGunhee21} and symmetrize polydisks of arbitrary complex dimensions \cite{GunheeChoYuanYuan20,NikolovNikolaiPflugPeterThomasPascalJZwonekWlodzimierz08} and also others \cite{GunheeCho20}. The proof of Theorem~\ref{thm:main} will be presented in Section 3.

Theorem~\ref{thm:main} holds on weakly pseudoconvex domains with very nice symmetries. In particular, if $\Omega$ is selected as in Theorem~\ref{th:CC}, and it would be particularly beneficial to obtain a distance comparison between the origin $(0,...,0)$ and the diagonal entry $(x,...,x)$. Especially, the following theorem is implicitly related to \cite[Problem 2.10]{JarnickiMarekPflugPeter13}. 

\begin{theo}\label{th:CC}
Let $n\geq2$. For each $0<\epsilon<\frac{1}{n^n}$, define
	\begin{equation*}
	\mathbb{D}^n_{\epsilon}:=\left\{ (z_1,...,z_n)\in{\mathbb{D}}^{n}:\left| \Pi^{n}_{i=1}z_i \right|<\epsilon \right\}.
	\end{equation*}	
	Then for any $(x,\cdots,x)\in \mathbb{D}^n_{\epsilon}$ satisfying $|x|>(\epsilon n)^{\frac{1}{n-1}}$, we have
	\[c_{\mathbb{D}^n_{\epsilon}}((0,...,0),(x,...,x)) <  d^{{KE}}_{\mathbb{D}^n_{\epsilon}}((0,...,0),(x,...,x)).\]	
\end{theo}

Note that $\epsilon>0$ must satisfy $\epsilon<\frac{1}{n^n}$ in order to take $(x,\cdots,x)\in \mathbb{D}^n_{\epsilon}$ with $|x|>(\epsilon n)^{\frac{1}{n-1}}$. The proof of Theorem~\ref{th:CC} will be presented in Section 4 with some additional remarks.

Lastly, one can easily extend the proof of Theorem~\ref{thm:main} with the complete K\"ahler--Einstein metric of Ricci curvature $-\lambda, \lambda>0$. 


\subsection*{Acknowledgements}
This work was partially supported by NSF grant DMS-1611745. I would like to thank my advisor Professor Damin Wu for many insights and helpful comments. I also would like to express my sincere thanks to the referees for reading the previous version of the paper carefully and for valuable comments to update the paper in a very progressive way with various detailed suggestions.

\section{Schwarz--Yau lemma}
The following generalized Schwarz lemma due to Yau will be used to compare the Carath\'eodory-distance and the K\"ahler--Einstein distance. 

\begin{theo}[the Schwarz-Yau lemma, \cite{YauShingTung78-1,RoydenHL80}]\label{th:SC}
	Let $(M,g)$ be a complete K\"ahler manifold with Ricci curvature bounded from below by a negative constant $K_1$. Let $(N,h)$ be another Hermitian manifold with holomorphic bisectional curvature bounded from above by a negative constant $K_2$. If there is a non-constant holomorphic map $f$ from $M$ to $N$, we have
	\[f^{*}h\leq \frac{K_1}{K_2}g. \] 
\end{theo}

We can use the upper bound of holomorphic sectional curvature of $(N,h)$ instead of the upper bound of bisectional curvature if $N$ is a Riemann surface \cite{RoydenHL80}. 

\section{$c_{\Omega}\lneq d^{{KE}}_{\Omega}$ for some complete Reinhardt domains $\Omega$ in $\mathbb{C}^n,n\geq 2$}

The proof of Theorem~\ref{thm:main} can be reduced to Lemma~\ref{lem:JP} which gives the comparison of the Carath\'eodory pseudo--distance and the K\"ahler--Einstein distance. Note that by Montel's theorem, for any two points in a complex manifold $M$, we can always achieve the extremal map $f\in \Hol(M, \mathbb{D})$ with respect to the Carath\'eodory pseudo--distance. 
In Lemma~\ref{lem:JP}, $\gamma_{\mathbb{D}}$ is the Poincar\'e metric on the unit disk and $\omega_{KE}(a)$ is the hermitian inner product on the holomorphic tangent space at $a\in \Omega$.

\begin{lem}\label{lem:JP}
	Let $\Omega$ be a bounded pseudoconvex domain in $\mathbb{C}^n$, and $a,b \in \Omega$, $a\neq b$. Suppose there exists $f\in \Hol(\Omega, \mathbb{D})$ that is extremal for $c_{\Omega}(a,b)$ such that 
	\begin{equation}\label{eq:assumption_of_Lemma}
	\gamma_{\mathbb{D}}(f(a),f'(a)X) <  \sqrt{(\omega_{KE}(a)(X,X))}, X \in \mathbb{C}^n-\left\{ 0\right\}. 
	\end{equation}	
	
	Then 
	\[c_{\Omega}(a,b) <  d^{{KE}}_{\Omega}(a,b).\]
	\begin{proof}
		By assumption, continuity of metrics gives $\epsilon, \delta$ such that 
		\[ \gamma_{\mathbb{D}}(f(z),f'(z)X) + \epsilon \| X\| \leq \sqrt{(\omega_{KE}(z)(X,X))}, z\in \mathbb{B}(a,\delta) \Subset \Omega, X\in \mathbb{C}^n.   \]
		On the other hand, by Theorem~\ref{th:SC}, we have 
		\[ \gamma_{\mathbb{D}} (f(z), f'(z)X) \leq \sqrt{(\omega_{KE}(z)(X,X))}, z\in \Omega, X\in \mathbb{C}^n.  \]	
		Let $\alpha : [0,1] \rightarrow \Omega$ be any piecewise $C^1$-curve joining $a$ and $b$. Denote by $t_0$ the maximal $t\in[0,1]$ such that $\alpha([0,t])\subset \overline{\mathbb{B}(a,\delta)}$. Denote the arc-length of $\alpha$ with respect to the K\"ahler--Einstein metric by $L_{KE}(\alpha)$. Then
		\[ L_{KE}(\alpha)= \int_{0}^{t_0}\sqrt{\omega_{KE}(\alpha(t))({\alpha'(t),\alpha'(t)})}dt + \int_{t_0}^{1}\sqrt{\omega_{KE}(\alpha(t))({\alpha'(t),\alpha'(t)})dt}  \] 	
		\[ \geq \int_{0}^{1}\gamma_{\mathbb{D}}(f(\alpha(t)), f'(\alpha(t))dt + \epsilon \int_{0}^{t_0} \| \alpha'(t)\| dt \]
		\[ \geq L_{\gamma_{\mathbb{D}}}(f\circ \alpha) +\epsilon \delta \geq \rho_{\mathbb{D}}(f(a),f(b)) + \epsilon \delta = c_{\Omega}(a,b)+\epsilon \delta.  \]
		Hence $d^{{KE}}_{\Omega}(a,b) \geq c_{\Omega}(a,b)+\epsilon \delta$.
		
	\end{proof}
\end{lem}

\begin{proof}[Proof of Theorem~\ref{thm:main}]
	
	The comparison between $c_{\Omega}$ and $d^{{KE}}_{\Omega}$ can be achieved once some extremal map $f\in \Hol(\Omega, \mathbb{D})$ with respect to the Carath\'eodory pseudo--distance satisfies the assumption of Lemma~\ref{lem:JP}. 
	
	Since $\Omega$ is a bounded domain in $\mathbb{C}^n$, we may assume that $\Omega$ is contained in the unit ball $B\subset \mathbb{C}^n$ centered at the origin due to the fact that the scaling transformation is a biholomorphism and the Carath\'eodory pseudo--distance and the K\"ahler--Einstein distance are preserved under any biholomorphism. 
	
	With the global coordinate $(z_1,...,z_n)\in \Omega$ in $\mathbb{C}^{n}$, let $\{\frac{\partial}{\partial{z_i}} | i=1,...,n\}$ be the basis on the holomorphic tangent bundle $T^{1,0}_{0}\Omega$. Without loss of generality, we may assume that $\{\frac{\partial}{\partial{z_i}}| i=1,...,n\}$ are orthonormal with repect to the Euclidean metric and orthogonal with respect to $\omega_{KE}(0)$. With the usual global coordinates $(z_1,\cdots,z_n)\in \Omega$, denote $\frac{\partial}{\partial{z_i}}=X_i,i=1,\cdots,n$. We may assume that $0<\omega_{KE}(0)(X_1,\overline X_1)\leq\omega_{KE}(0)(X_i,\overline X_i),i=1\cdots,n$. 
	
	We will show there exists a local hypersurface in $\Omega$ which is defined by  \[\sqrt{\omega_{KE}(X_1,\overline X_1)(0)}\sum^{n}_{k=1}z_k=1.\] From the assumption of Theorem~\ref{thm:main}, we can take $0<R\leq1$ such that $\Omega$ is contained in $B_R$ and the boundary of $\Omega$ touches the boundary of $B_R$. Then by using a projection map with respect to the first coordinate, we can define the holomorphic map from $(\Omega, \omega_{KE})$ to $(\mathbb{D}_R,h)$, where $h$ is the Poincar\'e metric on $\mathbb{D}_R$ the ball of radius $R$ in $\mathbb{C}^1$. Then by applying Theorem~\ref{th:SC}, $h(0)\leq R^2\omega_{KE}(0)$. This implies
	\begin{equation}\label{eq:fifth}
	1\leq \frac{1}{R} \leq R\sqrt{\omega_{KE}(0)(X_1,\overline X_1)}.
	\end{equation}
	In particular, $R \leq 1 \leq R\sqrt{\omega_{KE}(0)(X_1,\overline X_1)}$ and the assumption of the existence of $(p_1,\cdots,p_n)\in \Omega$ satisfying $\sum^n_{k=1}{p_k}=1$ implies the points $(z_1,\cdots,z_n)\in \Omega$ satisfying $\sqrt{\omega_{KE}(0)(X_1,\overline X_1)}\sum^{n}_{k=1}z_k=1$.  Here, the rescaling of the domain $\Omega\subset B$ and the fact that $\Omega$ is a complete Reinhardt domain are applied.
	
	We will control the extremal map with respect to the Carath\'eodory pseudo--distance $f : \Omega \rightarrow \mathbb{D}$ such that
	\begin{equation}\label{eq:third}
	c_{\Omega}((0,\cdots,0),(x,\cdots,x))=c_{\mathbb{D}}(f(0,\cdots,0),f(x,\cdots,x)).
	\end{equation} 
	After acting on the unit disk by an automorphism, we may assume that $f(0,\cdots,0)=0$. Also we can replace $f$ by the symmetrization map $\frac{1}{n!}\sum_{\sigma}f(z_{\sigma(1)},\cdots,z_{\sigma(n)})$ with all permutations $\sigma : \left\{1,\cdots,n\right\} \rightarrow \left\{1,\cdots,n\right\} $. Around the origin, one can write $f$ as a power series $\sum_{k=1}^{n} a_k z_k+\sum_{m=2}^{\infty}P_m(z_1,\cdots,z_n)=\sum_{i=1}^{n} a_i z_i+f_2(z_1,\cdots,z_n)$, where each $P_m(z_1,\cdots,z_n)$ is a homogeneous polynomial of degree $m$. Since we use $\frac{1}{n!}\sum_{\sigma}f(z_{\sigma(1)},\cdots,z_{\sigma(n)})$, we can replace $\sum_{i=1}^{n} a_i z_i$ by  $a \sum^{n}_{k=1}z_k$ for some real number $a$ (after acting on the unit disk by an automorphism if necessary). Then by Theorem~\ref{th:SC}, $f^* \gamma_{\mathbb{D}} \leq \sqrt{\omega_{KE}}$. In particular, we get
	\[a \leq \sqrt{\omega_{KE}(0)(X_1,\overline X_1)}.\]
	
	Now let's suppose $a=\sqrt{\omega_{KE}(0)(X_1,\overline X_1)}$. Since there exists a hypersurface in $\Omega$ given by $a\sum^{n}_{k=1}z_k=\sqrt{\omega_{KE}(X_1,\overline X_1)(0)}\sum^{n}_{k=1}z_k=1$, we can choose $(u_1,\cdots,u_n)\in \Omega$ such that $\sum^{n}_{k=1}u_k=\frac{1}{\sqrt{\omega_{KE}(0)(X_1,\overline X_1)}}$. Define $g : \mathbb{D} \rightarrow \Omega$ by $g(\lambda):=(\lambda u_1,\cdots,\lambda u_n)$. Since $\Omega$ is a complete Reinhardt domain, $g$ is a well-defined holomorphic function. Furtheremore, $f \circ g (0)=0$ and $(f \circ g)' (0)=1$. Thus $f \circ g = id_{\mathbb{D}}$ by the classical Schwarz's lemma that  $f_2(\lambda u_1, \cdots, \lambda u_n)=0$. Since $\lambda\in \mathbb{D}$ is arbitrary and with other choices of $(u_1,\cdots,u_n)$, we can take a small open set in $\Omega$ such that $f_2=0$ on this open set. Thus by the identity theorem, $f_2\equiv0$ on $\Omega$. Hence $f(z_1,\cdots,z_n)=a\sum^{n}_{k=1}z_k$, which is impossible since the image of $f$ can't hit the boundary of $\mathbb{D}$. Hence $a<\sqrt{\omega_{KE}(0)(X_1,\overline X_1)}$. In particular, this implies 
	\begin{equation}\label{eq:5}
	\left| f'(0)X_i \right| <  \sqrt{{\omega_{KE}(0)(X_i,\overline{X_i})}}, i=1,\cdots,n.
	\end{equation}
	Since we showed the above inequality for $X_i$'s, that are orthogonal to the Euclidean metric and $\omega_{KE}(0)$ at the same time, the same inequality holds for arbitrary tangent vectors $X \in \mathbb{C}^n-\left\{ 0\right\}$. Consequently, \eqref{eq:5} implies the assumption in Lemma~\ref{lem:JP}, and the proof is over.  
\end{proof}

\begin{rema}
	A separation between the Carath\'eodory pseudo--distance and the inner Carathe\'eodory pseudo--distance was made from similar argument (see \cite{VigueJeanPierre83}, \cite[Lemma 2.7.8]{JarnickiMarekPflugPeter13}). 
\end{rema}

\section{$c_{\mathbb{D}^n_{\epsilon}}((0,...,0),(x,...,x))\lneq d^{{KE}}_{\mathbb{D}^n_{\epsilon}}((0,...,0),(x,...,x))$}
Although the proof of the Thoerem~\ref{th:CC} contains the similar argument in the proof of the Theorem~\ref{thm:main}, we will provide the proof in detail, because the role of two fixed points $(x,...,x)$ and $(0,...,0)$ in $\mathbb{D}^n_{\epsilon}$ should be justified. 

\begin{proof}[Proof of Theorem~\ref{th:CC}]
	As we did in the proof of the Theorem~\ref{thm:main}, we establish the basic setting first. Fix $\epsilon>0$, with the global coordinate $(z_1,\cdots,z_n)\in \mathbb{D}^n_{\epsilon}$ in $\mathbb{C}^{n}$, let $\{\frac{\partial}{\partial{z_i}} | i=1,\cdots,z_n\}$ be the basis on the holomorphic tangent bundle $T^{1,0}_{0}\mathbb{D}^n_{\epsilon}$. Without loss of generality, we may assume that $\{\frac{\partial}{\partial{z_i}}| i=1,\cdots,n\}$ are orthonormal with repect to the Euclidean metric and orthogonal with respect to $\omega_{KE}(0)$. With the global coordinates $(z_1,\cdots,z_n)\in \mathbb{D}^n_{\epsilon}$, denote $\frac{\partial}{\partial{z_i}}=X_i, i=1,\cdots,n$. Also, we may assume that $0<\omega_{KE}(0)(X_1,\overline X_1)\leq\omega_{KE}(0)(X_i,\overline X_i),i=1,\cdots,n$. 
	For $z=(z_1,\cdots,z_n)$, define the holomorphic function $h : \mathbb{D}^n_{\epsilon} \rightarrow \mathbb{D}$ by $h(z):=\frac{1}{\epsilon}\Pi^{n}_{k=1}z_k$. 
	
	By the distance-decreasing property of the Carath\'eodory pseudo--distance,
\begin{equation*}
c_{\mathbb{D}}(h(0,\cdots,0),h(x,\cdots,x))\leq c_{\mathbb{D}^n_{\epsilon}}((0,\cdots,0),(x,\cdots,x)).
\end{equation*}	
Thus
	\begin{equation}\label{eq:first}
		c_{\mathbb{D}}(0,\frac{1}{\epsilon}{x}^{n})\leq c_{\mathbb{D}^n_{\epsilon}}((0,\cdots,0),(x,\cdots,x)).
	\end{equation}
	
From the hypothesis, take $(x,\cdots,x)\in \mathbb{D}^n_{\epsilon}$ satisfying $(\epsilon n)^{\frac{1}{n-1}}<|x|<1$.  	
We may assume $x$ is a positive real number so that $x \in \mathbb{D}$ satisfies
	\begin{equation}\label{eq:second}
		x > (\epsilon n)^{\frac{1}{n-1}}.
	\end{equation}
	
	Let $f : \mathbb{D}^n_{\epsilon} \rightarrow \mathbb{D}$ be the extremal map with respect to the Carath\'eodory pseudo--distance so that 
	\begin{equation}\label{eq:third}
		c_{\mathbb{D}^n_{\epsilon}}((0,\cdots,0),(x,\cdots,x))=c_{\mathbb{D}}(f(0,\cdots,0),f(x,\cdots,x)).
	\end{equation} 
We may assume that $f(0,\cdots,0)=0$. Since $\mathbb{D}^n_{\epsilon}$ is compatible with the symmetry by all permuations,	we may assume that $f$ itself is a symmetrization map so that $f(z)=a\sum_{k=1}^{n}z_k+\sum^{\infty}_{m=2} P_m(z)$, each $P_m$ is a homogeneous polynomial of degree $m$ and as in the same argument of the proof of Theorem~\ref{thm:main}, we get 
	\[a \leq \sqrt{\omega_{KE}(0)(X_1,\overline X_1)}.\]
	
	From the description of $\mathbb{D}^n_{\epsilon}$, we can find $(p_1,\cdots,p_n)\in \mathbb{D}^n_{\epsilon}$ satisfying $\sum^{n}_{k=1}z_k=1$ by taking one component with the magnitude almost one and the magnitudes of the other components almost zero. Then as in following the same argument of the proof of Theorem~\ref{thm:main}, there exists a local hypersurface of $\mathbb{D}^n_{\epsilon}$ given by \[\sqrt{\omega_{KE}(X_1,\overline X_1)(0)}\sum^{n}_{k=1}z_k=1.\] 
	Now we will show $a < \sqrt{\omega_{KE}(0)(X_1,\overline X_1)}$. Assume   $a=\sqrt{\omega_{KE}(0)(X_1,\overline X_1)}$. By \eqref{eq:first} and \eqref{eq:third},
	\[ c_{\mathbb{D}}(0,n\sqrt{\omega_{KE}(0)(X_1,\overline X_1)}x)=c_{\mathbb{D}^n_{\epsilon}}((0,0),(x,x))\geq c_{\mathbb{D}}(0,\frac{1}{\epsilon}{x}^{n}). \] 
	Thus we obtain $x \leq \left(n\epsilon \sqrt{(\omega_{KE}(0)(X_1,\overline{X_1})}\right)^{1/n}$. On the other hand, since the projection of $f : \mathbb{D}^n_{\epsilon} \rightarrow \mathbb{D}$ from $\mathbb{D}^n_{\epsilon}$ to the first coordinate induces the holomoprhic function from $\mathbb{D}$ to $\mathbb{D}$, the classical Schwarz lemma implies $\sqrt{(\omega_{KE}(0)(X_1,\overline{X_1})}=a\leq1$. In particular, we have
	\[n\epsilon\sqrt{\omega_{KE}(0)(X_1,\overline X_1)}=n\epsilon a\leq {n\epsilon}<1.   \]
	
	However, by \eqref{eq:second}, we also have $x > (\epsilon n)^{\frac{1}{n-1}}\geq \left(n\epsilon \sqrt{(\omega_{KE}(0)(X_1,\overline{X_1})}\right)^{1/n}$, which is impossible. Hence $a<\sqrt{\omega_{KE}(0)(X_1,\overline X_1)}$. 
	
	Then the rest of the proof follows as in the proof of Theorem~\ref{thm:main}.
\end{proof}

\bibliographystyle{spmpsci}
\bibliography{reference}

\end{document}